\documentclass[16pt, fleqn]{amsart}
\usepackage{amssymb,amsmath,amsthm,setspace,verbatim}
\usepackage{mathtools}
\usepackage{comment}
\usepackage[none]{hyphenat}
\onehalfspace

\newtheorem{thm}{Theorem}
\newtheorem*{thm*}{Sturm's Theorem}

\theoremstyle{definition}

\theoremstyle{remark}

\begin{document}

\title{Real Univariate Quintics}

\author{Elias Gonzalez}
\address{San Antonio, Texas}
\email{elias.gonzalez@nisd.net}

\author{David A. Weinberg}
\address{Department of Mathematics and Statistics, Texas Tech University, 
Lubbock, TX 79409}
\email{david.weinberg@ttu.edu}
\urladdr{www.math.ttu.edu/~dweinber}

\begin{abstract}
For the general monic quintic with real coefficients, polynomial conditions on the coefficients are derived as directly and as simply as possible from the Sturm sequence that will determine the real and complex root multiplicities together with the order of the real roots with respect to multiplicity.
\end{abstract}

\maketitle
\section{Quintics}
\subsection{Introduction.}
%%%%%%%%%%%%%%%%%%%%%%%%%%%%%%%%%%%%%%%%%%%%%%%%%%%%%%%%%%%%%%%%%%%%%%
\quad Consider the partition of the space of monic quintics, $x^{5}+px^{4}+qx^{3}+rx^{2}+sx+t$, according to the following root configurations:
\begin{enumerate}
\item 5 distinct real roots
\item 3 distinct real roots and 2 distinct complex conjugate roots
\item 1 real root and 4 distinct complex roots
\item 1 double real root and 3 distinct real roots
\begin{enumerate}
\item single root $<$ double root $<$ single root $<$ single root
\item double root $<$ single root $<$ single root $<$ single root
\item single root $<$ single root $<$ single root $<$ double root
\item single root $<$ single root $<$ double root $<$ single root
\end{enumerate}
\item 1 double real root and 1 single real root and 2 distinct complex conjugate roots
\begin{enumerate}
\item single root $<$ double root
\item double root $<$ single root
\end{enumerate}
\item 2 real double roots and 1 single real root
\begin{enumerate}
\item single root $<$ double root $<$ double root
\item double root $<$ single root $<$ double root
\item double root $<$ double root $<$ single root
\end{enumerate}
\item 2 complex conjugate double roots and 1 single real root
\item 1 triple root and 2 single real roots
\begin{enumerate}
\item triple root $<$ single root $<$ single root
\item single root $<$ triple root $<$ single root
\item single root $<$ single root $<$ triple root
\end{enumerate}
\item 1 triple root and 2 complex conjugate roots
\item 1 quadruple root and 1 single root
\begin{enumerate}
\item quadruple root $<$ single root
\item single root $<$ quadruple root
\end{enumerate}
\item 1 triple root and 1 double root
\begin{enumerate}
\item triple root $<$ double root
\item double root $<$ triple root
\end{enumerate}
\item 1 quintuple root
\end{enumerate}
We will find polynomial conditions on the coefficients $p$, $q$, $r$, $s$, and $t$ that will determine to which of these classes the polynomial belongs.  We will also determine polynomial conditions on the coefficients that will determine the order of the real roots with respect to multiplicity.  (Thus, there will be three sets of results: two for real polynomials and one for complex polynomials.)  In a previous paper [3], the authors gave the analogous results for cubics and quartics.

\quad If $c=\{c_{1},c_{2},$\ldots,$c_{m}\}$ is a finite sequence of real numbers, then the \emph{number of variations in sign} of $c$ is defined to be the number of $i$, $1\leq i\leq m-1$ such that $c_{i}c_{i+1}<0$, after dropping the 0's in $c$.  The \emph{Sturm sequence} for $f(x)$ is defined to be $f_{0}(x)=f(x)$, $f_{1}(x)=f'(x)$, $f_{2}(x)$, \ldots, $f_{s}(x)$, where\\
\begin{eqnarray*}
f_{0}(x)&=&q_{1}(x)f_{1}(x)-f_{2}(x), \quad \mathrm{deg} f_{2}(x)<\mathrm{deg} f_{1}(x) \\
& \vdots& \\
f_{i-1}(x)&=&q_{i}(x)f_{i}(x)-f_{i+1}(x), \quad \mathrm{deg} f_{i+1}(x)<\mathrm{deg} f_{i}(x)\\
& \vdots& \\
f_{s-1}(x)&=&q_{s}(x)f_{s}(x), \quad (f_{s+1}(x)=0)
\end{eqnarray*}

In other words, perform the Euclidean algorithm and change the sign of the remainder at each stage.  It should be noted here that there exists a more general definition of Sturm sequence, but the one given here suits our purpose.\\

\begin{thm*} Let $f(x)$ be a polynomial of positive degree in $\mathbb{R}[x]$ and let $f_{0}(x)=f(x)$, $f_{1}(x)=f'(x)$, \ldots, $f_{s}(x)$ be the Sturm sequence for $f(x)$.  Assume $[a,b]$ is an interval such that $f(a)\neq0$ and $f(b)\neq0$.  Then the number of distinct real roots of $f(x)$ in $(a,b)$ is $V_{a}-V_{b}$, where $V_{c}$ denotes the number of variations of sign of $\{f_{0}(c)$, $f_{1}(c)$, \ldots, $f_{s}(c)\}$.\end{thm*}

Proofs can be found in [4] and [6].  We can get the total number of real roots by looking at the limits as $a\to-\infty$ and $b\to+\infty$.  Thus, the total number of distinct real roots will depend only on the leading terms of the polynomials in the Sturm sequence.\\

Of course, in addition to Sturm's Theorem, we will also use the theorem that the complex roots of a real polynomial occur in conjugate pairs and the following (which the reader can prove as an exercise, or see [1], p. 102 (Cor. 4.1.9) or [2], p. 17 (Exercise 8)), less well known
\begin{thm} If the discriminant of a real polynomial is not zero, then the sign of the discriminant is $(-1)^{r}$, where $r$ is the number of complex conjugate pairs of roots.\end{thm}

In the 1990's, the authors Lu Yang, Songxian Liang, Jingzhong Zhang, Xiarong Hou, and Zhenbing Zeng studied root multiplicities and obtained extensive results.  In [5], [8], and [9], those authors developed and applied the notions of complete discrimination system, multiple factor sequence, and revised sign list, and proved general theorems about conditions for root multiplicities.

The significance of this paper is 1) the conditions giving the order of the real roots with respect to multiplicity are new and 2) the proof of the conditions for the real and complex root multiplicities is the simplest and most transparent possible.

\quad We use the computer software Maple to calculate the Sturm sequence for the quintic.  One obtains remainder polynomials of degrees 3, 2, 1, and 0, which we suggestively name $gcddeg3$, $gcddeg2$, $gcddeg1$, and $gcddeg0$, respectively.\\

\quad The Sturm sequence for the quintic is:

\begin{center}
$x^{5}+px^{4}+qx^{3}+rx^{2}+sx+t$\\
$5x^{4}+4px^{3}+3qx^{2}+2rx+s$\\
$gcddeg3=- \left( 2/5\,q-{\frac {4}{25}}\,{p}^{2} \right) {x}^{3}- \left( 3/5\,r-{\frac {3}{25}}\,pq \right) {x}^{2}- \left( 4/5\,s-{\frac {2}{25}}\,pr \right) x-t+1/25\,ps$\\

\tiny
\begin{eqnarray*}
gcddeg2 &=& -{\frac {25}{4}}\,{\frac { \left( -3\,{p}^{2}{q}^{2}+45\,{r}^{2}+8\,r{p}^{3}-40\,qs-38\,rpq+16\,{p}^{2}s+12\,{q}^{3} \right) {x}^{2}}{ \left( 2\,{p}^{2}-5\,q \right) ^{2}}} \\
\qquad \\
&-&{\frac {25}{2}}\,{\frac { \left( -25\,tq+4\,{q}^{2}r-3\,p{r}^{2}-21\,spq-{p}^{2}rq+10\,{p}^{2}t+6\,s{p}^{3}+30\,sr \right) x}{ \left( 2\,{p}^{2}-5\,q \right) ^{2}}} \\
\qquad \\
&-& {\frac{25}{4}}\,{\frac{-55\,tpq+16\,t{p}^{3}+75\,tr-3\,psr+4\,{q}^{2}s-{p}^{2}sq}{ \left( 2\,{p}^{2}-5\,q \right) ^{2}}}
\end{eqnarray*}
\normalsize
%%%%%%%%%%%%%%%%%%%%%%%%%%%%%%%%%
%$gcddeg2=-{\frac {25}{4}}\,{\frac { \left( -3\,{p}^{2}{q}^{2}+45\,{r}^{2}+8\,r{p}^{3}-40\,qs-38\,rpq+16\,{p}^{2}s+12\,{q}^{3} \right) {x}^{2}}{ \left( 2\,{p}^{2}-5\,q \right) ^{2}}}-{\frac {25}{2}}\,{\frac { \left( -25\,tq+4\,{q}^{2}r-3\,p{r}^{2}-21\,spq-{p}^{2}rq+10\,{p}^{2}t+6\,s{p}^{3}+30\,sr \right) x}{ \left( 2\,{p}^{2}-5\,q \right) ^{2}}}-{\frac{25}{4}}\,{\frac{-55\,tpq+16\,t{p}^{3}+75\,tr-3\,psr+4\,{q}^{2}s-{p}^{2}sq}{ \left( 2\,{p}^{2}-5\,q \right) ^{2}}}$\\
%%%%%%%%%%%%%%%%%%%%%%%%%%%%%%%%%%%%%%

\tiny
\[
\mathrm{gcddeg1} = \frac {8}{25}
\frac{\begin{aligned}
&( 6000\,{p}^{2}trqs-7500\,t{q}^{2}sr-
2490\,s{p}^{2}{q}^{2}{r}^{2}-528\,s{p}^{5}{q}^{2}r+1590\,s{p}^{3}{q}^{
3}r+2640\,{s}^{2}r{p}^{3}q+588\,s{r}^{2}{p}^{4}q-3300\,{q}^{2}{s}^{2}r
p \\
&-1550\,srp{q}^{4}-2400\,{p}^{3}tq{r}^{2}-584\,{p}^{6}tqr-3300\,{p}^{3
}t{q}^{2}s+1560\,{p}^{4}t{q}^{2}r+1080\,{p}^{5}tqs-1200\,{p}^{4}trs-
850\,{p}^{2}tr{q}^{3}\\
&+56\,{p}^{7}sqr+3000\,tp{q}^{2}{r}^{2}+3250\,tp{q
}^{3}s+4000\,{s}^{3}{q}^{2}-2200\,{s}^{2}{q}^{4}+3125\,{t}^{2}{q}^{3}-
1000\,tr{q}^{4}+2925\,s{r}^{2}{q}^{3}\\
&+108\,s{p}^{4}{q}^{4}+4185\,{s}^{
2}{p}^{2}{q}^{3}-2742\,{s}^{2}{p}^{4}{q}^{2}-315\,s{p}^{2}{q}^{5}-24\,
s{r}^{2}{p}^{6}-528\,{s}^{2}r{p}^{5}-3200\,q{s}^{3}{p}^{2}-460\,{p}^{3
}{r}^{3}{q}^{2}\\
&+4\,{p}^{6}{r}^{2}{q}^{2}-36\,{p}^{4}{r}^{2}{q}^{3}+540
\,{p}^{2}{r}^{4}q+450\,p{r}^{3}{q}^{3}+152\,{p}^{5}{r}^{3}q+105\,{p}^{
2}{r}^{2}{q}^{4}+216\,{p}^{5}t{q}^{3}-630\,{p}^{3}t{q}^{4}\\
&-24\,{p}^{7}
t{q}^{2}+480\,{p}^{5}t{r}^{2}+64\,{p}^{8}tr-112\,{p}^{7}ts-12\,{p}^{6}
s{q}^{3}+748\,{p}^{6}{s}^{2}q+600\,tp{q}^{5}-3750\,{t}^{2}{q}^{2}{p}^{
2}+1500\,{p}^{4}{t}^{2}q\\
&+640\,{p}^{4}{s}^{3}+300\,s{q}^{6}-108\,{p}^{4
}{r}^{4}-16\,{p}^{7}{r}^{3}-100\,{r}^{2}{q}^{5}-675\,{r}^{4}{q}^{2}-
200\,{p}^{6}{t}^{2}-72\,{p}^{8}{s}^{2}) x \end{aligned}
}{(-3\,{p}^{2}{
q}^{2}+45\,{r}^{2}+8\,r{p}^{3}-40\,qs-38\,rpq+16\,{p}^{2}s+12\,{q}^{3}
 )^{2}}
 \]
 \bigskip
 \[
 -\frac{4}{25} \frac{\begin{aligned}
 &(-5800\,tr{p}^{3}qs+7250\,t{q}
^{2}srp-1200\,{s}^{2}r{q}^{3}+3800\,t{p}^{2}{q}^{2}{r}^{2}+1224\,t{p}^
{5}{q}^{2}r-11850\,t{p}^{2}{q}^{3}s-3720\,t{p}^{3}{q}^{3}r+7560\,t{p}^
{4}{q}^{2}s\\
&-700\,t{r}^{2}{p}^{4}q+1160\,tr{p}^{5}s+8000\,tq{s}^{2}{p}^
{2}+3650\,trp{q}^{4}-332\,{p}^{4}{s}^{2}rq+1135\,{p}^{2}{s}^{2}{q}^{2}
r-450\,ps{r}^{2}{q}^{3}+460\,{p}^{3}s{q}^{2}{r}^{2}\\
&-4\,{p}^{6}s{q}^{2}
r+36\,{p}^{4}s{q}^{3}r-540\,{p}^{2}s{r}^{3}q-152\,{p}^{5}s{r}^{2}q-105
\,{p}^{2}sr{q}^{4}-2024\,{p}^{6}tsq-128\,{p}^{7}qrt-7500\,{p}^{2}{t}^{
2}qr\\
&-900\,t{q}^{6}+64\,{p}^{5}{s}^{3}+320\,{p}^{7}{t}^{2}-324\,t{p}^{4
}{q}^{4}+945\,t{p}^{2}{q}^{5}-4875\,t{r}^{2}{q}^{3}-16\,t{r}^{2}{p}^{6
}-10000\,t{q}^{2}{s}^{2}\\
&+6500\,t{q}^{4}s-1600\,t{p}^{4}{s}^{2}+108\,{p
}^{4}s{r}^{3}+16\,{p}^{7}s{r}^{2}+400\,p{s}^{3}{q}^{2}+300\,p{s}^{2}{q
}^{4}-315\,{p}^{3}{s}^{2}{q}^{3}+108\,{p}^{5}{s}^{2}{q}^{2}\\
&+28\,{p}^{6
}{s}^{2}r-320\,{p}^{3}{s}^{3}q+100\,r{q}^{5}s+675\,{r}^{3}{q}^{2}s-
6875\,{t}^{2}{q}^{3}p+7500\,{t}^{2}{q}^{2}{p}^{3}+9375\,{t}^{2}{q}^{2}
r+36\,{p}^{6}{q}^{3}t\\
&-2700\,{p}^{5}{t}^{2}q+1500\,{p}^{4}{t}^{2}r+192
\,{p}^{8}st-12\,{p}^{7}{s}^{2}q) \end{aligned}}{(-3\,{p}^{2}{q}^{2}+45\,{r}^{2
}+8\,r{p}^{3}-40\,qs-38\,rpq+16\,{p}^{2}s+12\,{q}^{3})^{2}}
\]
\normalsize

\bigskip

\tiny
\[
\mathrm{gcddeg0} =
\frac {25}{16} \frac{\begin{aligned} 
&( -3\,{p}^{2}{q}^{2}+45\,{r}^{2}+8\,r{
p}^{3}-40\,qs-38\,rpq+16\,{p}^{2}s+12\,{q}^{3} )^{2} ( 256
\,{p}^{5}{t}^{3}\\
&-192\,{p}^{4}qs{t}^{2}-128\,{p}^{4}{r}^{2}{t}^{2}+144
\,{p}^{4}r{s}^{2}t+144\,{p}^{3}{q}^{2}r{t}^{2}-6\,{p}^{3}{q}^{2}{s}^{2
}t-80\,{p}^{3}q{r}^{2}st+16\,{p}^{3}{r}^{4}t-27\,{p}^{2}{q}^{4}{t}^{2} \\
&+18\,{p}^{2}{q}^{3}rst-4\,{p}^{2}{q}^{2}{r}^{3}t-1600\,{p}^{3}q{t}^{3}
+160\,{p}^{3}rs{t}^{2}-36\,{p}^{3}{s}^{3}t+1020\,{p}^{2}{q}^{2}s{t}^{2
}+560\,{p}^{2}q{r}^{2}{t}^{2}-746\,{p}^{2}qr{s}^{2}t\\
&+24\,{p}^{2}{r}^{3
}st-630\,p{q}^{3}r{t}^{2}+24\,p{q}^{3}{s}^{2}t+356\,p{q}^{2}{r}^{2}st-
72\,pq{r}^{4}t+108\,{q}^{5}{t}^{2}-72\,{q}^{4}rst+16\,{q}^{3}{r}^{3}t\\
&+2000\,{p}^{2}r{t}^{3}-50\,{p}^{2}{s}^{2}{t}^{2}+2250\,p{q}^{2}{t}^{3}-
2050\,pqrs{t}^{2}+160\,pq{s}^{3}t-900\,p{r}^{3}{t}^{2}+1020\,p{r}^{2}{
s}^{2}t-900\,{q}^{3}s{t}^{2} \\
&+825\,{q}^{2}{r}^{2}{t}^{2}+560\,{q}^{2}r{
s}^{2}t-630\,q{r}^{3}st-2500\,ps{t}^{3}-3750\,qr{t}^{3}+2000\,q{s}^{2}
{t}^{2}+108\,{r}^{5}t-27\,{r}^{4}{s}^{2}+2250\,{r}^{2}s{t}^{2}\\
&-1600\,r
{s}^{3}t+256\,{s}^{5}+3125\,{t}^{4}-27\,{p}^{4}{s}^{4}+18\,{p}^{3}qr{s
}^{3}-4\,{p}^{3}{r}^{3}{s}^{2}-4\,{p}^{2}{q}^{3}{s}^{3}+{p}^{2}{q}^{2}
{r}^{2}{s}^{2}\\
&+144\,{p}^{2}q{s}^{4}-6\,{p}^{2}{r}^{2}{s}^{3}-80\,p{q}^
{2}r{s}^{3}+18\,pq{r}^{3}{s}^{2}+16\,{q}^{4}{s}^{3}-4\,{q}^{3}{r}^{2}{
s}^{2}-192\,pr{s}^{4}\\
&-128\,{q}^{2}{s}^{4}+144\,q{r}^{2}{s}^{3}
 ) \end{aligned}}
 {\begin{aligned} &( 14\,r{p}^{3}qs-62\,{q}^{2}srp-132\,p{s}^{2}r-28\,{p
}^{3}st-4\,{r}^{3}{p}^{3}-4\,{r}^{2}{q}^{3}-88\,{q}^{2}{s}^{2}+12\,{q}
^{4}s-18\,{p}^{4}{s}^{2}\\
&+125\,{t}^{2}q+{p}^{2}{q}^{2}{r}^{2}-3\,{p}^{2
}{q}^{3}s+117\,{r}^{2}qs+18\,{r}^{3}pq-6\,{r}^{2}{p}^{2}s+97\,q{s}^{2}
{p}^{2}-40\,t{q}^{2}r \\
&-300\,trs-27\,{r}^{4}+160\,{s}^{3}-66\,{p}^{2}qtr
+130\,ptsq-6\,{p}^{3}t{q}^{2}+120\,pt{r}^{2}+24\,p{q}^{3}t \\
&-50\,{p}^{2}
{t}^{2}+16\,{p}^{4}tr )^{2} (2\,{p}^{2}-5\,q )^{2} \end{aligned}
}
\]
\normalsize
\end{center}

Now let us do a Sturm analysis for this sequence.

The discriminant of the quintic is $D=256\,{p}^{5}{t}^{3}-192\,{p}^{4}qs{t}^{2}-128\,{p}^{4}{r}^{2}{t}^{2}+144\,{p}^{4}r{s}^{2}t+144\,{p}^{3}{q}^{2}r{t}^{2}-6\,{p}^{3}{q}^{2}{s}
^{2}t-80\,{p}^{3}q{r}^{2}st+16\,{p}^{3}{r}^{4}t-27\,{p}^{2}{q}^{4}{t}^{2}+18\,{p}^{2}{q}^{3}rst-4\,{p}^{2}{q}^{2}{r}^{3}t-1600\,{p}^{3}q{t}^
{3}+160\,{p}^{3}rs{t}^{2}-36\,{p}^{3}{s}^{3}t+1020\,{p}^{2}{q}^{2}s{t}^{2}+560\,{p}^{2}q{r}^{2}{t}^{2}-746\,{p}^{2}qr{s}^{2}t+24\,{p}^{2}{r}
^{3}st-630\,p{q}^{3}r{t}^{2}+24\,p{q}^{3}{s}^{2}t+356\,p{q}^{2}{r}^{2}st-72\,pq{r}^{4}t+108\,{q}^{5}{t}^{2}-72\,{q}^{4}rst+16\,{q}^{3}{r}^{3
}t+2000\,{p}^{2}r{t}^{3}-50\,{p}^{2}{s}^{2}{t}^{2}+2250\,p{q}^{2}{t}^{3}-2050\,pqrs{t}^{2}+160\,pq{s}^{3}t-900\,p{r}^{3}{t}^{2}+1020\,p{r}^{
2}{s}^{2}t-900\,{q}^{3}s{t}^{2}+825\,{q}^{2}{r}^{2}{t}^{2}+560\,{q}^{2}r{s}^{2}t-630\,q{r}^{3}st-2500\,ps{t}^{3}-3750\,qr{t}^{3}+2000\,q{s}^
{2}{t}^{2}+108\,{r}^{5}t-27\,{r}^{4}{s}^{2}+2250\,{r}^{2}s{t}^{2}-1600\,r{s}^{3}t+256\,{s}^{5}+3125\,{t}^{4}-27\,{p}^{4}{s}^{4}+18\,{p}^{3}q
r{s}^{3}-4\,{p}^{3}{r}^{3}{s}^{2}-4\,{p}^{2}{q}^{3}{s}^{3}+{p}^{2}{q}^{2}{r}^{2}{s}^{2}+144\,{p}^{2}q{s}^{4}-6\,{p}^{2}{r}^{2}{s}^{3}-80\,p{
q}^{2}r{s}^{3}+18\,pq{r}^{3}{s}^{2}+16\,{q}^{4}{s}^{3}-4\,{q}^{3}{r}^{2}{s}^{2}-192\,pr{s}^{4}-128\,{q}^{2}{s}^{4}+144\,q{r}^{2}{s}^{3}$.\\

Let $L3$, $L2$, and $L1$ denote the appropriate factors of the leading coefficients of gcddeg3, gcddeg2, and gcddeg1, respectively:\\
\begin{eqnarray*}
L3&=&2p^{2}-5q \\
L2&=&3\,{p}^{2}{q}^{2}-45\,{r}^{2}-8\,r{p}^{3}+40\,qs+38\,rpq-16\,{p}^{2}s-12\,{q}^{3} \\
L1&=&-54\,{r}^{4}+320\,{s}^{3}+250\,q{t}^{2}-176\,{q}^{2}{s}^{2}+24\,{q}^{4
}s-36\,{p}^{4}{s}^{2}-100\,{p}^{2}{t}^{2} \\
&-&124\,srp{q}^{2}+28\,sr{p}^{3
}q+260\,sptq-132\,{p}^{2}qrt-12\,{p}^{3}t{q}^{2}+48\,pt{q}^{3} \\
&+&2\,{p}^
{2}{q}^{2}{r}^{2}-6\,{q}^{3}s{p}^{2}-80\,{q}^{2}rt+194\,q{s}^{2}{p}^{2
}-600\,str-264\,p{s}^{2}r+36\,{r}^{3}pq \\
&-&56\,s{p}^{3}t-12\,s{r}^{2}{p}^
{2}+234\,sq{r}^{2}+32\,{p}^{4}tr+240\,p{r}^{2}t-8\,{q}^{3}{r}^{2}-8\,{
r}^{3}{p}^{3}
\end{eqnarray*}

Notice that the numerator of $gcddeg0$ is $25(L2)^{2}D$.

\subsection{Real and Complex Root Multiplicities.}
The number of real roots is determined by the leading coefficients of the Sturm sequence.  However, in constructing a proof, one has to be careful about whether or not some of these leading coefficients are zero.\\

\subsubsection{5 Distinct Roots.}
%%%%%%%%Exhibit Sturm analysis.%%%%%%%%%%%%%%%

The Sturm analysis shows that if $L3>0$, $L2>0$, $L1>0$, and $D>0$, then there are 5 distinct real roots.  (Clearly, there is no other way to get 5 real roots.)  By Theorem 1, the condition for 3 real single roots and two complex conjugate roots is $D<0$.  Again by Theorem 1, the conditions for one real root and two distinct pairs of complex conjugate roots are $D>0$ and $(L3\leq0$ or $L2\leq0$ or $L1\leq0)$.\\

\subsubsection{1 Double and 3 Single Roots.}
We have $D=0$ and $L1\neq0$.  There are 2 or 4 distinct real roots.  First consider the case where $L2\neq0$ and $L3\neq0$.  An immediate and superficial consequence of Sturm's Theorem is that to get 4 real roots, one must have $L3>0$, $L2>0$, and $L1>0$.  Otherwise there are 2 real roots.  But there is an interesting point here.  A careful analysis of the Sturm sequence (as shown below) shows that it suffices to consider $L1$ alone!  Here is the Sturm analysis in this case.\\

\begin{center}
\begin{tabular}{c c c c c c}
 & $x^{5}$ & $x^{4}$ & $(L3)x^{3}$ & $(L2)x^{2}$ & $(L1)x$\\
\hline
$-\infty$ & $-$ & $+$ & $-L3$ & $L2$ & $-L1$\\
$\infty$ & $+$ & $+$ & $L3$ & $L2$ & $L1$\\
\end{tabular}
\end{center}

\begin{center}
\begin{tabular}{c c c c c l}
$L3$ & $L2$ & $L1$ & & & \\
\hline
$+$ & $+$ & $+$ & & $4 - 0$ & $ = 4 $\\
$+$ & $+$ & $-$ & & $3 - 1$ & $ = 2 $\\
$+$ & $-$ & $+$ & & $2 - 2$ & $ = 0 $ (Impossible)\\
$+$ & $-$ & $-$ & & $3 - 1$ & $ = 2 $\\
$-$ & $+$ & $+$ & & $2 - 2$ & $ = 0 $ (Impossible)\\
$-$ & $+$ & $-$ & & $1 - 3$ & $ = -2 $ (Impossible)\\
$-$ & $-$ & $+$ & & $2 - 2$ & $ = 0 $ (Impossible)\\
$-$ & $-$ & $-$ & & $3 - 1$ & $ = 2 $\\
\end{tabular}
\end{center}

Thus, the Sturm analysis shows that there are 4 distinct real roots when $L1>0$, and 2 distinct real roots when $L1<0$.  In the special cases where $L3=0$ or $L2=0$, it is impossible to get four sign changes in the Sturm analysis table, hence, there must be only 2 distinct real roots.  We would like to verify that the conditions in these special cases still depend only on the original $L1$ (defined above).\\

Consider the special case where $L2=0$ and $L3\neq0$.  We know (by the superficial consequence of Sturm's Theorem) that there can only be 2 real roots.  But we would like to verify that this case forces $L1<0$.  (Note that here we use $\widetilde{L1}$, as it is not the same coefficient of $x$ as in the original Sturm sequence due to the special case $L2=0$.)\\ 

\begin{center}
\begin{tabular}{c c c c c}
 & $x^{5}$ & $x^{4}$ & $(L3)x^{3}$ & $(\widetilde{L1})x$\\
\hline
$-\infty$ & $-$ & $+$ & $-L3$ & $-\widetilde{L1}$\\
$\infty$ & $+$ & $+$ & $L3$ & $\widetilde{L1}$\\
\end{tabular}
\end{center}

\begin{center}
\begin{tabular}{c c c c l}
$L3$ & $\widetilde{L1}$ & & & \\
\hline
$+$ & $+$ & & $2 - 0$ & $ = 2 $\\
$+$ & $-$ & & $3 - 1$ & $ = 2 $\\
$-$ & $+$ & & $1 - 1$ & $ = 0 $ (Impossible)\\
$-$ & $-$ & & $2 - 2$ & $ = 0 $ (Impossible)\\
\end{tabular}
\end{center}

The Sturm analysis shows that $L3$ has to be greater than zero.  If $L2=0$, then $s=-1/8\,{\frac {-3\,{p}^{2}{q}^{2}+8\,{p}^{3}r-38\,rpq+45\,{r}^{2}+12\,{q}^{3}}{2\,{p}^{2}-5\,q}}$ and substituting this into $L1$ gives 
\begin{equation*}
-{\frac{1}{64}}\,{\frac{ \splitfrac{(520\,{q}^{3}r+135\,{p}^{3}{q}^{3}-18\,{p}^{5}{q}^{2}+800\,{p}^{2}qt-252\,p{q}^{4}-1000\,t{q}^{2}-160\,{p}^{4}t}{+1350\,{r}^{3}+604\,{p}^{2}{q}^{2}r-380\,{p}^{4}qr-2205\,pq{r}^{2}+558\,{p}^{3}{r}^{2}+48\,{p}^{6}r )^{2}}}{ ( 2\,{p}^{2}-5\,q )^{3}}}.
\end{equation*}  
(Notice the perfect square in the numerator, a result of computer algebra.)  Therefore, $L1<0$.\\

Now consider the special case where $L3=0$ and $L2\neq0$.  Let us recompute the Sturm sequence for this case (solve $L3=0$ for $q$ and substitute $q=2p^2/5$ into the original quintic).
\begin{center}
${x}^{5}+p{x}^{4}+2/5\,{p}^{2}{x}^{3}+r{x}^{2}+sx+t$\\
$5\,{x}^{4}+4\,p{x}^{3}+6/5\,{p}^{2}{x}^{2}+2\,rx+s$\\
$new1gcddeg2 = - \left( -{\frac {6}{125}}\,{p}^{3}+3/5\,r \right) {x}^{2}- \left( -{
\frac {2}{25}}\,pr+4/5\,s \right) x+1/25\,ps-t$\\
$new1gcddeg1 = 
-{\dfrac {50}{27}}\,
{\dfrac { \left( -2700\,sr{p}^{5}-82500\,p{s}^{2}r-
187500\,str+15000\,s{p}^{3}t+25500\,s{r}^{2}{p}^{2}+75000\,p{r}^{2}t \atop -
10500\,{p}^{4}tr+4200\,{p}^{4}{s}^{2}+72\,{p}^{8}s+360\,{p}^{7}t-60\,{
p}^{6}{r}^{2}-16875\,{r}^{4}+100000\,{s}^{3}+2000\,{r}^{3}{p}^{3}
 \right) x}{ \left( 2\,{p}^{3}-25\,r \right) ^{3}}}-{\dfrac {25}{27}}\,
{\dfrac { \left( -60\,{p}^{6}sr-16875\,s{r}^{3}+7625\,{p}^{2}{s}^{2}r+51250\,t{
p}^{2}{r}^{2}+10500\,{p}^{4}st-6600\,t{p}^{5}r \atop +18750\,{p}^{3}{t}^{2}-
234375\,r{t}^{2}+2000\,{p}^{3}s{r}^{2}-181250\,psrt+216\,t{p}^{8}+
250000\,t{s}^{2}-10000\,p{s}^{3}-450\,{p}^{5}{s}^{2} \right)}
{ \left( 2\,{p}^{
3}-25\,r \right) ^{3}}}$
\end{center}
\small
\begin{align*}
&new1gcddeg0 =  \bigg( -{\dfrac {27}{500}}\, ( 2\,{p}^{3}-25\,r )^{3}
 (161250\,{p}^{2}{s}^{3}{r}^{2}+7031250\,s{r}^{2}{t}^{2}-5000000
\,{s}^{3}rt-84375\,{s}^{2}{r}^{4}\\
&+337500\,t{r}^{5}+800000\,{s}^{5}+
9765625\,{t}^{4}  -2160\,{p}^{8}str+78000\,{p}^{5}{r}^{2}st-202500\,{p}^
{4}{s}^{2}tr \\
&-2062500\,{p}^{3}s{t}^{2}r-712500\,{p}^{2}st{r}^{3}+
3187500\,p{s}^{2}{r}^{2}t+1296\,{t}^{2}{p}^{10}+480\,{p}^{8}{s}^{3}-
75000\,{p}^{5}{t}^{3} \\
&+31625\,{p}^{4}{s}^{4}-54000\,{p}^{7}r{t}^{2}+
1800\,{p}^{7}t{s}^{2}+90000\,{p}^{6}{t}^{2}s+1200\,{p}^{6}{r}^{3}t-300
\,{p}^{6}{r}^{2}{s}^{2} \\
&-17500\,{p}^{5}r{s}^{3}+712500\,{p}^{4}{r}^{2}{
t}^{2}-40000\,{p}^{3}t{r}^{4}+10000\,{p}^{3}{s}^{2}{r}^{3}+87500\,{p}^
{3}t{s}^{3}+1562500\,{p}^{2}{t}^{3}r\\
&+2343750\,{p}^{2}{t}^{2}{s}^{2}-
2812500\,p{t}^{2}{r}^{3}-600000\,p{s}^{4}r-7812500\,p{t}^{3}s \bigg) \bigg/ \bigg( -2700\,sr{p}^{5} \\
&-82500\,p{s}^{2}r-187500\,str+15000\,s{p}^{3
}t+25500\,s{r}^{2}{p}^{2}+75000\,p{r}^{2}t-10500\,{p}^{4}tr \\
&+4200\,{p}^
{4}{s}^{2}+72\,{p}^{8}s+360\,{p}^{7}t-60\,{p}^{6}{r}^{2}-16875\,{r}^{4
}+100000\,{s}^{3}+2000\,{r}^{3}{p}^{3}\bigg)^{2}
\end{align*}
\normalsize

\begin{center}
\begin{tabular}{c c c c c}
 & $x^{5}$ & $x^{4}$ & $(\widetilde{L2})x^{2}$ & $(\widetilde{L1})(\widetilde{L2})x$\\
\hline
$-\infty$ & $-$ & $+$ & $\widetilde{L2}$ & $-(\widetilde{L1})(\widetilde{L2})$\\
$\infty$ & $+$ & $+$ & $\widetilde{L2}$ & $(\widetilde{L1})(\widetilde{L2})$\\
\end{tabular}
\end{center}

\begin{center}
\begin{tabular}{c c c c c l}
$\widetilde{L2}$ & $\widetilde{L1}$ & $(\widetilde{L1})(\widetilde{L2})$ & & & \\
\hline
$+$ & $+$ & $+$ & & $2 - 0$ & $ = 2 $\\
$+$ & $-$ & $-$ & & $1 - 1$ & $ = 0 $ (Impossible)\\
$-$ & $+$ & $-$ & & $3 - 1$ & $ = 2 $\\
$-$ & $-$ & $+$ & & $2 - 2$ & $ = 0 $ (Impossible)\\
\end{tabular}
\end{center}

$\widetilde{L1}$ and $\widetilde{L2}$ denote the appropriate factors of the leading coefficients in $new1gcddeg1$ and $new1gcddeg2$.  It follows from the Sturm analysis that we must have $\widetilde{L1} > 0$.  Substitute $q=2p^2/5$ into $L1$: the result is $L1=-{\frac {216}{25}}\,sr{p}^{5}-264\,p{s}^{2}r-600\,str+48\,s{p}^{3}t+{
\frac {408}{5}}\,s{r}^{2}{p}^{2}+240\,p{r}^{2}t-{\frac {168}{5}}\,{p}^
{4}tr+{\frac {32}{5}}\,{r}^{3}{p}^{3}+{\frac {336}{25}}\,{p}^{4}{s}^{2
}+{\frac {144}{625}}\,{p}^{8}s+{\frac {144}{125}}\,{p}^{7}t-{\frac {24
}{125}}\,{p}^{6}{r}^{2}-54\,{r}^{4}+320\,{s}^{3}$.
  We see that $\widetilde{L1}$ and $L1$ have opposite signs.  This shows that the condition is $L1 < 0$.\\
  Finally consider the case $L2=0$ and $L3=0$.  Recomputing the Sturm sequence and doing the Sturm analysis yields $L1^*=p^4-125s > 0$, where $L1^*$ is the appropriate coefficient of the resulting $gcddeg1$.  Furthermore, computer algebra shows that $L3=0$ and $L2=0$ imply that the original $L1=-\frac{64(p^4-125s)^3}{390625}$.  Thus, we must have $L1<0$.  Therefore, in all cases, the condition for one double root, one real root, and two complex conjugate roots is $L1<0$.

\subsubsection{One Triple and Two Single Roots or Two Double and One Single Root.}
We have $D=0$ and $L1=0$ and $L2\neq0$. (By definition of Sturm sequence and discriminant, it follows as a consequence that the constant term of $gcddeg1$ must be zero.)\\
First consider the generic case $L3\neq0$.  The discriminant of gcddeg2 will tell whether there is a triple root or 2 double roots, and in the latter case whether they are real or complex conjugate.
Here is the appropriate factor (having the same sign as) of the discriminant of gcddeg2:\\
$D2 = 160\,{s}^{2}{q}^{3}+{p}^{4}{q}^{2}{r}^{2}-236\,{p}^{4}q{s}^{2}-136\,{p
}^{5}st-3\,{p}^{4}{q}^{3}s+48\,{p}^{5}{q}^{2}t-1100\,t{q}^{3}r-128\,{p
}^{6}rt-45\,sp{r}^{3}+60\,s{q}^{2}{r}^{2}-1380\,{p}^{3}{r}^{2}t-8\,{p}
^{2}{q}^{3}{r}^{2}+337\,{p}^{2}{q}^{2}{s}^{2}+408\,{p}^{3}r{s}^{2}-12
\,{p}^{4}{r}^{2}s-357\,{p}^{3}{q}^{3}t+6\,{p}^{3}q{r}^{3}+24\,{p}^{2}{
q}^{4}s-500\,{p}^{2}q{t}^{2}-24\,p{q}^{2}{r}^{3}+660\,p{q}^{4}t+36\,{p
}^{6}{s}^{2}+100\,{p}^{4}{t}^{2}+9\,{p}^{2}{r}^{4}-48\,s{q}^{5}+625\,{
t}^{2}{q}^{2}+1500\,tqsr-4\,{p}^{5}qrs+1028\,{p}^{4}qrt+11\,{p}^{3}{q}
^{2}rs+800\,{p}^{3}qst-1735\,{p}^{2}{q}^{2}rt+20\,p{q}^{3}rs+5475\,pq{
r}^{2}t-600\,s{p}^{2}tr-1150\,stp{q}^{2}-1380\,{s}^{2}rpq-3\,s{p}^{2}{
r}^{2}q-3375\,{r}^{3}t+16\,{r}^{2}{q}^{4}+900\,{r}^{2}{s}^{2}$.\\

If $D2=0$, then there is a triple root and 2 single roots.\\
If $D2>0$, then there are 2 double real roots and one single real root.\\
If $D2<0$, then there are 2 complex conjugate double roots and one single real root.\\

When there is a triple root, it is necessary to do a Sturm analysis.

\begin{center}
\begin{tabular}{c c c c c}
 & $x^5$ & $x^4$ & $(L3)x^3$ & $(L2)x^2$\\
\hline
$-\infty$ & $-$ & $+$ & $-L3$ & $L2$\\
$\infty$ & $+$ & $+$ & $L3$ & $L2$\\
\end{tabular}
\end{center}

\begin{center}
\begin{tabular}{c c c c l}
$L3$ & $L2$ & & & \\
\hline
$+$ & $+$ & & $3 - 0$ & $ = 3 $\\
$+$ & $-$ & & $2 - 1$ & $ = 1 $\\
$-$ & $+$ & & $1 - 2$ & $ = -1 $ (Impossible)\\
$-$ & $-$ & & $2 - 1$ & $ = 1 $ \\
\end{tabular}
\end{center}

Thus, if $L2>0$, then there are 3 real roots, and if $L2<0$, then there is 1 real root.\\

Now consider the special case $L3=0$.\\

Then by the Sturm Theorem, there can only be 1 real root (one cannot get 3 sign changes).

%%%%%%%DISCUSS GCDDEG2's DISCRIM%%%%%%%%%%%%%%
(Only the triple root case needs the Sturm analysis because if there are 2 double roots, then the sign of $D2$ tells whether they are real or complex.)  
%%%%%%%%%%%%%EXHIBIT STURM ANALYSIS%%%%%%%%%%%

If $L3=0$, solve for $q$ and recompute the Sturm sequence: $q=2p^2/5$.  The result is\\

${x}^{5}+p{x}^{4}+2/5\,{p}^{2}{x}^{3}+r{x}^{2}+sx+t$\\

$5\,{x}^{4}+4\,p{x}^{3}+6/5\,{p}^{2}{x}^{2}+2\,rx+s$\\

$- \left( -{\frac {6}{125}}\,{p}^{3}+3/5\,r \right) {x}^{2}- \left( -{
\frac {2}{25}}\,pr+4/5\,s \right) x+1/25\,ps-t$\\

$-{\dfrac{50}{27}}\,{\dfrac{ \left( 15000\,s{p}^{3}t-2700\,{p}^{5}rs-
187500\,srt-82500\,pr{s}^{2}-10500\,{p}^{4}rt+25500\,{p}^{2}{r}^{2}s+
75000\,p{r}^{2}t\atop +2000\,{p}^{3}{r}^{3}-60\,{p}^{6}{r}^{2}+100000\,{s}^{
3}+360\,{p}^{7}t+72\,s{p}^{8}-16875\,{r}^{4}+4200\,{s}^{2}{p}^{4}
 \right) x}{ \left( 2\,{p}^{3}-25\,r \right) ^{3}}} \\ -{\dfrac {25}{27}}\, 
{\frac{ \left( 7625\,{p}^{2}{s}^{2}r-16875\,{r}^{3}s+10500\,{p}^{4}st+51250\,
t{p}^{2}{r}^{2}-6600\,t{p}^{5}r+2000\,{p}^{3}s{r}^{2}-60\,{p}^{6}sr- \atop
181250\,psrt  -10000\,p{s}^{3}+216\,t{p}^{8}-234375\,r{t}^{2}+250000\,t{
s}^{2}+18750\,{p}^{3}{t}^{2}-450\,{p}^{5}{s}^{2} \right)}{ ( 2\,{p}^{3}-
25\,r )^3}}$

Let us denote by $\widetilde{L2}$ and $\widetilde{L1}$ the appropriate factors of the leading coefficients of $gcddeg2$ and $gcddeg1$ from the new Sturm sequence.  We need $\widetilde{L2}\neq0$ and $\widetilde{L1}=0$.  If we plug in $q=2p^2/5$ into the appropriate factor of the numerator of the leading coefficient of the original $gcddeg2$, then we get $-{\frac {9}{125}}\, \left( 2\,{p}^{3}-25\,r \right) ^{2}$.  If we plug in $q=2p^2/5$ into the appropriate factor of the numerator of the leading coefficient of the original $gcddeg1$, then we get\\
$-300\,srt-132\,pr{s}^{2}+{\frac {204}{5}}\,{p}^{2}{r}^{2}s+120\,p{r}^{
2}t+24\,s{p}^{3}t-{\frac {108}{25}}\,{p}^{5}rs-{\frac {84}{5}}\,{p}^{4
}rt+160\,{s}^{3}+{\frac {16}{5}}\,{p}^{3}{r}^{3}+{\frac {72}{125}}\,{p
}^{7}t+{\frac {72}{625}}\,s{p}^{8}+{\frac {168}{25}}\,{s}^{2}{p}^{4}-{
\frac {12}{125}}\,{p}^{6}{r}^{2}-27\,{r}^{4}$.\\

The Sturm analysis shows that if $L2>0$, there are 3 real roots, and if $L2<0$ there is 1 real root.\\

Notice that in the special case that $L3=0$, $L2$ is forced to be negative.\\

In the special case $L3=0$, it is also necessary to compare the discriminant of the new gcddeg2 with that of the original gcddeg2.  The discriminant of the new gcddeg2 is\\
$-{\frac {24}{3125}}\,s{p}^{4}+{\frac {24}{125}}\,{p}^{3}t+{\frac {4}{
625}}\,{p}^{2}{r}^{2}-{\frac {4}{125}}\,prs+{\frac {16}{25}}\,{s}^{2}-
{\frac {12}{5}}\,rt$.\\

If we plug $q=2p^2/5$ into the appropriate factor ($D2$) of the original gcddeg2, we get\\
$-{\frac {36}{5}}\, \left( 2\,{p}^{3}-25\,r \right) ^{2} \left( 6\,s{p}
^{4}-150\,{p}^{3}t-5\,{p}^{2}{r}^{2}+25\,prs+1875\,rt-500\,{s}^{2}
 \right)$.\\
 
Thus, $D2$ still distinguishes between the triple root and two double roots, even in the special case $L3=0$.

%%%%%%%WORKSHEET STUFF%%%%%%%%%%%%%%%%%%%%%%%

\subsubsection{1 Quadruple and 1 Single or 1 Triple and 1 Double}
We have $D=0$, $L1=0$, $L2=0$, and $L3\neq0$.  In this case, all roots are real.\\
In this case, we know that gcddeg3 must have a multiple root (either triple or double).  Here is the Sturm sequence for $gcddeg3$:\\

\begin{align*}
& gcddeg3 = - \left( 2/5\,q-{\frac {4}{25}}\,{p}^{2} \right) {x}^{3}- \left( 3/5\,r-{\frac {3}{25}}\,pq \right) {x}^{2}- \left( 4/5\,s-{\frac {2}{25}}\,pr \right) x-t+1/25\,ps\\
& (gcddeg3)' = -3\, \left( 2/5\,q-{\frac {4}{25}}\,{p}^{2} \right) {x}^{2}-2\, \left( 3/5\,r-{\frac {3}{25}}\,pq \right) x-4/5\,s+{\frac {2}{25}}\,pr
\end{align*}
\begin{align*}
& 3gcddeg1 = -{\dfrac {1}{75}}\,{\dfrac { \left( 8\,{p}^{3}r-80\,s{p}^{2}-3\,{p}^{2}{q}^{2}+10\,prq+200\,sq-75\,{r}^{2} \right) x}{-5\,q+2\,{p}^{2}}} \\
&-{\dfrac{1}{75}}\,{\dfrac{6\,{p}^{3}s-150\,{p}^{2}t-{p}^{2}rq-5\,qps+5\,p{r}^{2}+375\,qt-50\,sr}{-5\,q+2\,{p}^{2}}}
\end{align*}
\begin{align*}
& 3gcddeg0 = -{\dfrac {2}{25}}\, \left( -5\,q+2\,{p}^{2} \right)  \bigg( -
5400\,s{p}^{5}t-108\,{p}^{5}srq+2700\,{p}^{4}rqt+45\,{p}^{3}sr{q}^{2} \\
&+ 135000\,{p}^{2}srt-32000\,{s}^{3}{p}^{2}+108\,{s}^{2}{p}^{6}+67500\,{p
}^{4}{t}^{2}-337500\,{p}^{2}q{t}^{2} \\
&+ 33750\,ps{q}^{2}t-16875\,p{r}^{2}
qt-337500\,srqt+1275\,{p}^{2}sq{r}^{2}+3375\,{p}^{2}r{q}^{2}t\\
&-1500\,pr
{s}^{2}q+80000\,{s}^{3}q-22500\,{s}^{2}{r}^{2}+84375\,{r}^{3}t+421875
\,{q}^{2}{t}^{2}\\
&+10\,{p}^{3}{r}^{3}q-420\,{p}^{4}s{r}^{2}+540\,{p}^{4}
{s}^{2}q+27\,{p}^{4}{q}^{3}s-9\,{p}^{4}{r}^{2}{q}^{2}-675\,{p}^{3}{q}^
{3}t \\
&-13500\,{p}^{3}{r}^{2}t +4200\,{p}^{3}r{s}^{2}-2925\,{p}^{2}{s}^{2}
{q}^{2}+1125\,ps{r}^{3}+32\,{p}^{5}{r}^{3}-225\,{p}^{2}{r}^{4}
 \bigg) \Bigg/  \\
 & \left( 8\,{p}^{3}r-80\,s{p}^{2}-3\,{p}^{2}{q}^{2}+10\,prq+
200\,sq-75\,{r}^{2} \right)^{2}\\
\end{align*}

Therefore, by the results for the cubic [3], the conditions for a quadruple root are $D3={\frac {864}{15625}}\,s{p}^{5}t+{\frac {432}{390625}}\,{p}^{5}srq-{\frac {432}{15625}}\,{p}^{4}rqt-{\frac {36}{78125}}\,{p}^{3}sr{q}^{2}-{\frac {864}{625}}\,{p}^{2}srt+{\frac {1024}{3125}}\,{s}^{3}{p}^{2}-{
\frac {432}{390625}}\,{s}^{2}{p}^{6}-{\frac {432}{625}}\,{p}^{4}{t}^{2}+{\frac {432}{125}}\,{p}^{2}q{t}^{2}-{\frac {216}{625}}\,ps{q}^{2}t+{
\frac {108}{625}}\,p{r}^{2}qt+{\frac {432}{125}}\,srqt-{\frac {204}{15625}}\,{p}^{2}sq{r}^{2}-{\frac {108}{3125}}\,{p}^{2}r{q}^{2}t+{
\frac {48}{3125}}\,pr{s}^{2}q-{\frac {512}{625}}\,{s}^{3}q+{\frac {144}{625}}\,{s}^{2}{r}^{2}-{\frac {108}{125}}\,{r}^{3}t-{\frac {108}{25}}
\,{q}^{2}{t}^{2}-{\frac {8}{78125}}\,{p}^{3}{r}^{3}q+{\frac {336}{78125}}\,{p}^{4}s{r}^{2}-{\frac {432}{78125}}\,{p}^{4}{s}^{2}q-{\frac 
{108}{390625}}\,{p}^{4}{q}^{3}s+{\frac {36}{390625}}\,{p}^{4}{r}^{2}{q}^{2}+{\frac {108}{15625}}\,{p}^{3}{q}^{3}t+{\frac {432}{3125}}\,{p}^{
3}{r}^{2}t-{\frac {672}{15625}}\,{p}^{3}r{s}^{2}+{\frac {468}{15625}}\,{p}^{2}{s}^{2}{q}^{2}-{\frac {36}{3125}}\,ps{r}^{3}-{\frac {128}{
390625}}\,{p}^{5}{r}^{3}+{\frac {36}{15625}}\,{p}^{2}{r}^{4}=0$ (the discriminant of $gcddeg3$), $M1=\left( 8\,{p}^{3}r-80\,s{p}^{2}-3\,{p}^{2}{q}^{2}+10\,prq+200\,sq-75\,{r}^{2} \right) =0$, where $M1$ is the leading coefficient of $3gcddeg1$ above.  The condition for triple and double root are $D3=0$ and $M1\neq0$. (But it is unnecessary to say $D3=0$ because we know gcddeg3 MUST have a multiple root.)\\

\subsubsection{1 Quintuple Root}
The conditions are $D=0$, $L1=0$, $L2=0$, and $L3=0$.\\

\subsection{Order of the Real Roots with Respect to Multiplicity.}

\subsubsection{1 Double and 3 Single Roots}
We begin with the case where we have 1 double real root and 3 single real roots.  Solving $gcddeg1$ for the double root yields:\\

\begin{center}
$doubleroot=\frac{C_{0}}{C_{1}}$\\
where $C_{0}=48\,s{p}^{4}t+4\,s{p}^{3}{r}^{2}+80\,{p
}^{3}{t}^{2}-32\,{p}^{3}rqt-3\,{p}^{3}{s}^{2}q+7\,{s}^{2}{p}^{2}r-{p}^
{2}sr{q}^{2}-4\,{p}^{2}{r}^{2}t+9\,{p}^{2}t{q}^{3}-266\,sq{p}^{2}t+16
\,p{s}^{3}+146\,ptr{q}^{2}-18\,sp{r}^{2}q+290\,sptr-275\,pq{t}^{2}+12
\,p{s}^{2}{q}^{2}+4\,s{q}^{3}r-195\,{r}^{2}qt+260\,s{q}^{2}t+27\,s{r}^
{3}+375\,{t}^{2}r-36\,{q}^{4}t-48\,r{s}^{2}q-400\,t{s}^{2}$\\
and $C_{1}=2\,(-28\,s{p}^{3}t-50\,{p}^{2}{t}^{2}+120\,p{r}^{2}t-132\,{s}^{2}rp-6\,{p}^{3}{q}^{2
}t+14\,{p}^{3}rsq-62\,pr{q}^{2}s-4\,{p}^{3}{r}^{3}-18\,{s}^{2}{p}^{4}-
88\,{s}^{2}{q}^{2}+12\,s{q}^{4}-4\,{r}^{2}{q}^{3}+125\,q{t}^{2}+97\,{s
}^{2}{p}^{2}q-6\,s{p}^{2}{r}^{2}-3\,{p}^{2}{q}^{3}s+{p}^{2}{q}^{2}{r}^
{2}+18\,p{r}^{3}q+117\,sq{r}^{2}-40\,r{q}^{2}t-300\,srt-27\,{r}^{4}+
160\,{s}^{3}+130\,qpst-66\,{p}^{2}rqt+24\,p{q}^{3}t+16\,{p}^{4}tr)$\\
come from the coefficients of $gcddeg1$.\\
\end{center}
%%%%Note that C_1 is the SAME as L1, so it is not needed.%%%%
We now translate the double root to the origin by substituting $x=y+doubleroot$ into the quintic.  The resulting quintic then must have a factor of $y^{2}$, as we have shifted the double root to the origin.  Here is the remaining cubic factor:\\

\begin{center}
$leftovercubic = {y}^{3}+ \left( p+5\,{\frac {C_{{0}}}{C_{{1}}}} \right) {y}^{2}+
 \left( q+4\,{\frac {pC_{{0}}}{C_{{1}}}}+10\,{\frac {{C_{{0}}}^{2}}{{C
_{{1}}}^{2}}} \right) y+3\,{\frac {qC_{{0}}}{C_{{1}}}}+r+6\,{\frac {p{
C_{{0}}}^{2}}{{C_{{1}}}^{2}}}+10\,{\frac {{C_{{0}}}^{3}}{{C_{{1}}}^{3}
}}$\\
\end{center}

\quad We now need information about how many real roots of a cubic $x^3 + px^2 +qx +r$ are positive. If we perform a Sturm sequence for the cubic, but this time checking the variations in sign from 0 to $\infty$, we will find out how many of these single roots are positive.  The signs at 0 and $\infty$ are determined by the constant terms and by the leading coefficients.  The cases actually condense, and depend only on the constant terms in the Sturm sequence.\\
\quad It actually takes several pages to write down the proof, using Sturm's Theorem, but the result is that if we are in the case of three single real roots, then if $q>0$, $r>0$, and $pq-9r>0$, then there are 0 positive single roots, and if $q>0$, $r<0$, and $pq-9r<0$, then there are 3 positive single roots.  If neither of these cases hold, then, if $r>0$, then there are 2 positive single roots, and if $r<0$, then there is 1 positive single root.  If we are in the case of one single real root, then if $r>0$, then there are 0 positive single roots, and if $r<0$, there is 1 positive single root.\\

\quad By these results, when there are three single real roots, then if $\left( q+4\,{\frac {pC_{{0}}}{C_{{1}}}}+10\,{\frac {{C_{{0}}}^{2}}{{C_{{1}}}^{2}}} \right)>0 $,  $3\,{\frac {qC_{{0}}}{C_{{1}}}}+r+6\,{\frac {p{C_{{0}}}^{2}}{{C_{{1}}}^{2}}}+10\,{\frac {{C_{{0}}}^{3}}{{C_{{1}}}^{3}}} > 0$, and ${\frac {p{C_{{1}}}^{3}q+4\,{p}^{2}{C_{{1}}}^{2}C_{{0}}-24\,p{C_{{0}}}^
{2}C_{{1}}-22\,qC_{{0}}{C_{{1}}}^{2}-40\,{C_{{0}}}^{3}-9\,r{C_{{1}}}^{
3}}{{C_{{1}}}^{3}}}
 > 0$ , then there are 0 positive roots; if $\left(q+4\,{\frac{pC_{{0}}}{C_{{1}}}}+10\,{\frac {{C_{{0}}}^{2}}{{C_{{1}}}^{2}}} \right) > 0$, $3\,{\frac {qC_{{0}}}{C_{{1}}}}+r+6\,{\frac {p{C_{{0}}}^{2}}{{C_{{1}}}^{2}}}+10\,{\frac {{C_{{0}}}^{3}}{{C_{{1}}}^{3}}} < 0$, and ${\frac {p{C_{{1}}}^{3}q+4\,{p}^{2}{C_{{1}}}^{2}C_{{0}}-24\,p{C_{{0}}}^
{2}C_{{1}}-22\,qC_{{0}}{C_{{1}}}^{2}-40\,{C_{{0}}}^{3}-9\,r{C_{{1}}}^{
3}}{{C_{{1}}}^{3}}}
 < 0$, then there are 3 positive roots.  If neither of these cases hold, then if $3\,{\frac {qC_{{0}}}{C_{{1}}}}+r+6\,{\frac {p{C_{{0}}}^{2}}{{C_{{1}}}^{2}}}+10\,{\frac {{C_{{0}}}^{3}}{{C_{{1}}}^{3}}} > 0$, then there are 2 positive roots, and if $3\,{\frac {qC_{{0}}}{C_{{1}}}}+r+6\,{\frac {p{C_{{0}}}^{2}}{{C_{{1}}}^{2}}}+10\,{\frac {{C_{{0}}}^{3}}{{C_{{1}}}^{3}}} < 0$, then there is 1 positive root.\\

\subsubsection{1 Double and 1 Single}
$doubleroot$ and $leftovercubic$ are the same as in the preceding section.  Then by the results on the number of positive or negative roots for a cubic, if the constant term of $leftovercubic$ is positive, then the double root of the quintic is bigger than the single root of the quintic, while if the constant term of $leftovercubic$ is negative, then the single root of the quintic is bigger than the double root of the quintic, i.e.\\
if $3\,{\frac {qC_{{0}}}{C_{{1}}}}+r+6\,{\frac {p{C_{{0}}}^{2}}{{C_{{1}}}^{2}}}+10\,{\frac {{C_{{0}}}^{3}}{{C_{{1}}}^{3}}}>0$, then double root $>$ single root, while\\
if $3\,{\frac {qC_{{0}}}{C_{{1}}}}+r+6\,{\frac {p{C_{{0}}}^{2}}{{C_{{1}}}^{2}}}+10\,{\frac {{C_{{0}}}^{3}}{{C_{{1}}}^{3}}}<0$, then single root $>$ double root.\\
But now we must consider the special cases where $L2=0$ or $L3=0$.  In each of these special cases, it is necessary to consider the recomputed Sturm sequences.  Computer algebra shows that in each of these special cases, the calculation of the double root and the leftovercubic from the recomputed Sturm sequences coincide with the result of substituting for $s$, respectively for $q$, (solving for s in $L2=0$, respectively for $q$ in $L3=0$), in $doubleroot$ and $leftovercubic$ from the original Sturm sequence.  Therefore, the conditions above for the relative position of the double and single root hold in all cases. 

\subsubsection{1 Triple and 2 Single}
In this case the greatest common divisor of the quintic and its derivative has degree 2.  First consider the generic case where $L3\neq0$.  Recall that the discriminant of $gcddeg2$ has been denoted $D2$.  If $D2=0$, then there is a triple root.  This triple root is $-(1/2)$ times the coefficient of $x$ in the monic associate of gcddeg2.\\

\begin{center}
$tripleroot={\frac{C_{2,1}}{L2}}$\\
\end{center}
where $C_{2,1}=6\,s{p}^{3}+4\,{q}^{2}r-3\,p{r}^{2}-21\,spq+30\,sr+10\,{p}^{2}t-25\,tq-{p}^{2}rq$\\
Now translate the triple root to the origin via $y=x-tripleroot$.  The resulting quintic will have a factor of $y^3$.  By computer algebra, the remaining quadratic factor is\\
$leftoverquadratic={y}^{2}+ \left( p+5\,{\frac {C_{{2,1}}}{{\it L2}}} \right) y+q+4\,{\frac {pC_{{2,1}}}{{\it L2}}}+10\,{\frac {{C_{{2,1}}}^{2}}{{{\it L2}}^
{2}}}
$\\

The appropriate factor of the discriminant of this leftoverquadratic is\\
$D22=-4\,q{{L2}}^{2}-6\,pC_{{2,1}}{L2}-15\,{C_{{2,1}}}^{2}+{{L2
}}^{2}{p}^{2}
   $\\
If $D22 < 0$, then the two single roots are complex.  If $D22 > 0$, then the two single roots are real.  If the constant term of $leftoverquadratic$ is negative, then the configuration is single $<$ triple $<$ single.  If the constant term of $leftoverquadratic$ is positive, then if the coefficient of $y$ in $leftoverquadratic$ is negative, then triple $<$ single $<$ single, while if the coefficient of $y$ in $leftoverquadratic$ is positive, then single $<$ single $<$ triple.\\

Now consider the special case $L3=0$.  It is necessary to consider the recomputed Sturm sequence.  But it follows from Sturm's Theorem that it is impossible for the two single roots to be real.  It is interesting to note that computer algebra shows that the discriminant of gcddeg2 in the recomputed Sturm sequence must be negative.

\subsubsection{2 Double and 1 Single}
Again, in this case the greatest common divisor of the quintic and its derivative must have degree 2.  If $D2<0$, then the two double roots are complex conjugate.  If $D2>0$, then the two double roots are real.  $(gcddeg2)^2$ must divide the quintic.  We can get the single root by solving for x in the quotient of the quintic divided by $(gcddeg2)^2$.\\

\begin{center}
$singleroot={\frac {C_{3}}{L2}}$.\\
where $C_{3}=-34\,{p}^{2}rq+8\,r{p}^{4}+44\,spq-3\,{p}^{3}{q}^{2}-8\,s{p}^
{3}+12\,p{q}^{3}+57\,p{r}^{2}-16\,{q}^{2}r+100\,tq-120\,sr-40\,{p}^{2}
t$\\
\end{center}

Then we can translate singleroot to the origin in $gcddeg2$ by substituting $x=y+singleroot$ into $gcddeg2$.  The relative positions of the 2 double roots and 1 single root can then be determined from the coefficients of the resulting quadratic polynomial.  This quadratic polynomial is the following:\\
$\frac {25}{4}\,(L2)\,y^2+ \left( -\frac {25}{2}\,\,C_{2, 1}+\frac {25}{2}\,\,C_{3} \right) y+\left(\frac {25}{4}\,\frac {C_{3}^{2}}{(L2)}+\frac {25}{4}\,C_{2, 0}-\frac {25}{2}\,\frac {C_{3}C_{2, 1}}{(L2)}\right)$\\
where $C_{2, 1} = -{p}^{2}rq+4\,{q}^{2}r-25\,tq-21\,spq+30\,sr+10\,{p}^{2}t-3\,p{r}^{2}+6\,s{p}^{3}$\\
and $C_{2, 0}= -16\,t{p}^{3}-75\,tr+3\,psr-4\,{q}^{2}s+55\,tpq+{p}^{2}sq$\\
come from the coefficients of $gcddeg2's$ linear and constant terms.  After making this quadratic polynomial monic, we can say that if the resulting constant term is negative, then we have double $<$ single $<$ double.  If the resulting constant term is positive, then if the resulting coeff of x is negative, we have single $<$ double $<$ double, and if the resulting coeff of x is positive, we have double $<$ double $<$ single.  Furthermore note that only the case L3 NOT equal to 0 is relevant because if it is zero, then by the Sturm analysis, you cannot have 3 real roots.
%%%%%  When writing the paper, be sure to identify all the factors in this quadratic in y - they all come from the original Sturm sequence or from singleroot. %%%%%%

\subsubsection{1 Quadruple and 1 Single}
In this case, the greatest common divisor of the quintic and its derivative has degree 3.  A quadruple root of the quintic is a triple root of the greatest common divisor.  The quadruple root is thus $-(1/3)$ times the coefficient of $x^2$ in the monic associate of gcddeg3.\\
$quadrupleroot=1/2\,{\frac {5\,r-pq}{-5\,q+2\,{p}^{2}}}$.\\
We can then translate the quadruple root to the origin in the original quintic and get a leftover linear factor whose root is the translated single root.\\
By computer algebra, we obtain that the translated single root is $-1/2\,{\frac {-15\,pq+4\,{p}^{3}+25\,r}{-5\,q+2\,{p}^{2}}}$.\\
Therefore, if $-1/2\,{\frac {-15\,pq+4\,{p}^{3}+25\,r}{-5\,q+2\,{p}^{2}}}>0$, then we have quadruple $<$ single, while if $-1/2\,{\frac {-15\,pq+4\,{p}^{3}+25\,r}{-5\,q+2\,{p}^{2}}}<0$, then we have single $<$ quadruple.

\subsubsection{1 Triple and 1 Double}
Again, the greatest common divisor of the quintic and its derivative has degree 3.  If this greatest common divisor does not have a triple root, then it must have a double root (because the original quintic cannot have 3 double roots), and the quintic must have 1 triple root and 1 double root.  We can now apply the results from our previous paper, [3], that determine the conditions on the cubic that will determine whether double $<$ single or single $<$ double, to gcddeg3, and this will give us the conditions on the quintic for triple $<$ double or double $<$ triple.  By computer algebra, the result is that if ${\dfrac {27}{4}}\,{\dfrac {\left(-1250\,t{q}^{2}+1000\,tq{p}^{2}-200\,t{p}^{4}
+50\,ps{q}^{2}-40\,{p}^{3}sq+8\,{p}^{5}s+375\,{r}^{2}q-150\,{r}^{2}{p}
^{2}\atop-150\,rp{q}^{2}+60\,r{p}^{3}q+15\,{p}^{2}{q}^{3}-6\,{p}^{4}{q}^{2}
-125\,{r}^{3}+75\,{r}^{2}pq-15\,r{p}^{2}{q}^{2}+{p}^{3}{q}^{3}\right)}{
 \left( -5\,q+2\,{p}^{2} \right) ^{3}}} > 0$, then we have double $>$ triple, while if 
 
 ${\dfrac {27}{4}}\,{\dfrac {\left(-1250\,t{q}^{2}+1000\,tq{p}^{2}-200\,t{p}^{4}
+50\,ps{q}^{2}-40\,{p}^{3}sq+8\,{p}^{5}s+375\,{r}^{2}q-150\,{r}^{2}{p}
^{2}\atop-150\,rp{q}^{2}+60\,r{p}^{3}q+15\,{p}^{2}{q}^{3}-6\,{p}^{4}{q}^{2}
-125\,{r}^{3}+75\,{r}^{2}pq-15\,r{p}^{2}{q}^{2}+{p}^{3}{q}^{3}\right)}{
 \left( -5\,q+2\,{p}^{2} \right) ^{3}}} < 0$, then we have triple $>$ double.\newpage
 
 \section{Summary}

\textbf{Notation.}\\
\begin{align*}
D&=256\,{p}^{5}{t}^{3}-192\,{p}^{4}qs{t}^{2}-128\,{p}^{4}{r}^{2}{t}^{2}+144\,{p}^{4}r{s}^{2}t+144\,{p}^{3}{q}^{2}r{t}^{2}-6\,{p}^{3}{q}^{2}{s}
^{2}t\\
&-80\,{p}^{3}q{r}^{2}st+16\,{p}^{3}{r}^{4}t-27\,{p}^{2}{q}^{4}{t}^{2}+18\,{p}^{2}{q}^{3}rst-4\,{p}^{2}{q}^{2}{r}^{3}t-1600\,{p}^{3}q{t}^
{3}\\
&+160\,{p}^{3}rs{t}^{2}-36\,{p}^{3}{s}^{3}t+1020\,{p}^{2}{q}^{2}s{t}^{2}+560\,{p}^{2}q{r}^{2}{t}^{2}-746\,{p}^{2}qr{s}^{2}t+24\,{p}^{2}{r}
^{3}st\\
&-630\,p{q}^{3}r{t}^{2}+24\,p{q}^{3}{s}^{2}t+356\,p{q}^{2}{r}^{2}st-72\,pq{r}^{4}t+108\,{q}^{5}{t}^{2}-72\,{q}^{4}rst+16\,{q}^{3}{r}^{3
}t\\
&+2000\,{p}^{2}r{t}^{3}-50\,{p}^{2}{s}^{2}{t}^{2}+2250\,p{q}^{2}{t}^{3}-2050\,pqrs{t}^{2}+160\,pq{s}^{3}t-900\,p{r}^{3}{t}^{2}\\
&+1020\,p{r}^{
2}{s}^{2}t-900\,{q}^{3}s{t}^{2}+825\,{q}^{2}{r}^{2}{t}^{2}+560\,{q}^{2}r{s}^{2}t-630\,q{r}^{3}st-2500\,ps{t}^{3}-3750\,qr{t}^{3}\\
&+2000\,q{s}^
{2}{t}^{2}+108\,{r}^{5}t-27\,{r}^{4}{s}^{2}+2250\,{r}^{2}s{t}^{2}-1600\,r{s}^{3}t+256\,{s}^{5}+3125\,{t}^{4}-27\,{p}^{4}{s}^{4}\\
&+18\,{p}^{3}q
r{s}^{3}-4\,{p}^{3}{r}^{3}{s}^{2}-4\,{p}^{2}{q}^{3}{s}^{3}+{p}^{2}{q}^{2}{r}^{2}{s}^{2}+144\,{p}^{2}q{s}^{4}-6\,{p}^{2}{r}^{2}{s}^{3}\\
&-80\,p{
q}^{2}r{s}^{3}+18\,pq{r}^{3}{s}^{2}+16\,{q}^{4}{s}^{3}-4\,{q}^{3}{r}^{2}{s}^{2}\\
&-192\,pr{s}^{4}-128\,{q}^{2}{s}^{4}+144\,q{r}^{2}{s}^{3}; \text{the discriminant of the quintic}\\
& \text{Let L3, L2, and L1 denote the appropriate factors of the leading coefficients of} \\ 
& \text{gcddeg3, gcddeg2, and gcddeg1, respectively:}\\
\end{align*}
\begin{align*}
L3&=2p^2-5q; \text {appropriate factor of the leading coefficient of gcddeg3}\\
L2&=40\,qs-16\,{p}^{2}s-8\,r{p}^{3}+38\,rpq+3\,{p}^{2}{q}^{2}-12
\,{q}^{3}-45\,{r}^{2}\\
L1&=-264\,p{s}^{2}r-12\,{p}^{3}t{q}^{2}+36\,{r}^{3}pq-124\,srp{q}
^{2}+28\,sr{p}^{3}q+260\,sptq-132\,{p}^{2}qrt\\
&+240\,p{r}^{2}t+234\,sq{r
}^{2}+32\,{p}^{4}tr+48\,pt{q}^{3}-56\,s{p}^{3}t-80\,{q}^{2}rt+194\,q{s
}^{2}{p}^{2}-600\,str\\
&-6\,{q}^{3}s{p}^{2}+2\,{p}^{2}{q}^{2}{r}^{2}-12\,
s{r}^{2}{p}^{2}-54\,{r}^{4}+320\,{s}^{3}-8\,{q}^{3}{r}^{2}-8\,{r}^{3}{
p}^{3}+250\,q{t}^{2}\\
&-176\,{q}^{2}{s}^{2}+24\,{q}^{4}s-36\,{p}^{4}{s}^{
2}-100\,{p}^{2}{t}^{2}\\
D2&=24\,{p}^{2}{q}^{4}s-1100\,{q}^{3}rt+800\,{p}^{3}qst-1735\,{p}
^{2}{q}^{2}rt-3\,{p}^{2}q{r}^{2}s+20\,p{q}^{3}rs-600\,{p}^{2}rst\\
&-1150
\,p{q}^{2}st+5475\,pq{r}^{2}t-1380\,pqr{s}^{2}+1500\,qrst+6\,{p}^{3}q{
r}^{3}+{p}^{4}{q}^{2}{r}^{2}-128\,{p}^{6}rt\\
&+660\,p{q}^{4}t-136\,{p}^{5
}st-3\,{p}^{4}{q}^{3}s-236\,{p}^{4}q{s}^{2}+337\,{p}^{2}{q}^{2}{s}^{2}
+48\,{p}^{5}{q}^{2}t-357\,{p}^{3}{q}^{3}t\\
&-12\,{p}^{4}{r}^{2}s-45\,p{r}
^{3}s+60\,{q}^{2}{r}^{2}s-8\,{p}^{2}{q}^{3}{r}^{2}-500\,{p}^{2}q{t}^{2
}-24\,p{q}^{2}{r}^{3}-1380\,{p}^{3}{r}^{2}t\\
&+408\,{p}^{3}r{s}^{2}-4\,{p
}^{5}qrs+1028\,{p}^{4}qrt+11\,{p}^{3}{q}^{2}rs+36\,{p}^{6}{s}^{2}+100
\,{p}^{4}{t}^{2}+9\,{p}^{2}{r}^{4}\\
&-48\,{q}^{5}s+16\,{q}^{4}{r}^{2}+160
\,{q}^{3}{s}^{2}+625\,{q}^{2}{t}^{2}-3375\,{r}^{3}t+900\,{r}^{2}{s}^{2};\\
& \text{the discriminant of gcddeg2}\\
\end{align*}
\begin{align*}
M1&= 8\,{p}^{3}r-80\,s{p}^{2}-3\,{p}^{2}{q}^{2}+10\,prq+200\,sq-75\,{r}^{2};\\
& \text{appropriate factor of the leading coefficient of 3gcddeg1} \\
& \text{(from the Sturm sequence for gcddeg3)}\\
C_{0}&=48\,s{p}^{4}t+4\,s{p}^{3}{r}^{2}+80\,{p
}^{3}{t}^{2}-32\,{p}^{3}rqt-3\,{p}^{3}{s}^{2}q+7\,{s}^{2}{p}^{2}r-{p}^
{2}sr{q}^{2}\\
&-4\,{p}^{2}{r}^{2}t+9\,{p}^{2}t{q}^{3}-266\,sq{p}^{2}t+16
\,p{s}^{3}+146\,ptr{q}^{2}-18\,sp{r}^{2}q+290\,sptr-275\,pq{t}^{2}\\
&+12
\,p{s}^{2}{q}^{2}+4\,s{q}^{3}r-195\,{r}^{2}qt+260\,s{q}^{2}t+27\,s{r}^
{3}\\
&+375\,{t}^{2}r-36\,{q}^{4}t-48\,r{s}^{2}q-400\,t{s}^{2};\\
& \text{appropriate factor of the constant term of gcddeg1.} \\
\end{align*}
\begin{align*}
C_{1}&=2\,(-28\,s{p}^{3}t-50\,{p}^{2}{t}^{2}+120\,p{r}^{2}t-132\,{s}^{2}rp-6\,{p}^{3}{q}^{2
}t+14\,{p}^{3}rsq-62\,pr{q}^{2}s\\
&-4\,{p}^{3}{r}^{3}-18\,{s}^{2}{p}^{4}-
88\,{s}^{2}{q}^{2}+12\,s{q}^{4}-4\,{r}^{2}{q}^{3}+125\,q{t}^{2}+97\,{s
}^{2}{p}^{2}q\\
&-6\,s{p}^{2}{r}^{2}-3\,{p}^{2}{q}^{3}s+{p}^{2}{q}^{2}{r}^
{2}+18\,p{r}^{3}q+117\,sq{r}^{2}-40\,r{q}^{2}t-300\,srt\\
&-27\,{r}^{4}+
160\,{s}^{3}+130\,qpst-66\,{p}^{2}rqt+24\,p{q}^{3}t+16\,{p}^{4}tr) \\
%%%%Note that C_1 is the SAME as L1 so it is not needed.%%%%
C_{2,1}&=6\,s{p}^{3}+4\,{q}^{2}r-3\,p{r}^{2}-21\,spq+30\,sr+10\,{p}^{2}t-25\,tq-{p}^{2}rq;\\
&\text{numerator of tripleroot} \\
D22&=-4\,q{{L2}}^{2}-6\,pC_{{2,1}}{L2}-15\,{C_{{2,1}}}^{2}+{{L2
}}^{2}{p}^{2};\\
& \text{appropriate factor of the discriminant of leftoverquadratic} \\
\end{align*}
\begin{align*}
C_{3}&=-34\,{p}^{2}rq+8\,r{p}^{4}+44\,spq-3\,{p}^{3}{q}^{2}-8\,s{p}^
{3}+12\,p{q}^{3}+57\,p{r}^{2}-16\,{q}^{2}r+100\,tq \\
& -120\,sr-40\,{p}^{2}t; \text{numerator of singleroot} \\
C_{2, 0}&= -16\,t{p}^{3}-75\,tr+3\,psr-4\,{q}^{2}s+55\,tpq+{p}^{2}sq;\\
& \text{appropriate factor of the constant term of gcddeg2} \\
C_{4}&=-1/2\,{\frac {-15\,pq+4\,{p}^{3}+25\,r}{-5\,q+2\,{p}^{2}}};\text{quadrupleroot}\\
\end{align*}
\begin{align*}
C_{5}&={\frac {27}{4}}\,{\frac { \left(-1250\,t{q}^{2}+1000\,tq{p}^{2}-200\,t{p}^{4}
+50\,ps{q}^{2}-40\,{p}^{3}sq+8\,{p}^{5}s+375\,{r}^{2}q-150\,{r}^{2}{p}
^{2}-150\,rp{q}^{2} \atop+60\,r{p}^{3}q +15\,{p}^{2}{q}^{3}-6\,{p}^{4}{q}^{2}
-125\,{r}^{3}+75\,{r}^{2}pq-15\,r{p}^{2}{q}^{2}+{p}^{3}{q}^{3}\right)}{
 \left( -5\,q+2\,{p}^{2} \right) ^{3}}}\\
 %%%%From Sturm analysis of gcddeg3.%%%%
\end{align*}
\begin{align*}
F_{1}&={q+\frac{4pC_{0}}{C_{1}}+10\frac{C_{0}^{2}}{C_{1}^2}}\\
F_{2}&={\frac{3qC_{0}}{C_{1}}+r+\frac{pC_{0}^{2}}{C_{1}^{2}}+\frac{10C_{0}^{3}}{C_{1}^{3}}}\\
F_{3}&={\frac{pC_{1}^{3}q+4p^{2}C_{1}^{2}C_{0}-24pC_{0}^{2}C_{1}-22qC_{0}C_{1}^{2}-40C_{0}^{3}-9rC_{1}^{3}}{C_{1}^{3}}}\\
\end{align*}
\begin{align*}
F_{4}&={\frac{C_{3}^{2}}{L2^{2}}+\frac{C_{2,0}}{L2}-2\frac{2C_{3}C_{2,1}}{L2^{2}}}\\
F_{5}&={\frac{-2C_{2,1}}{L2}+\frac{2C_{3}}{L2}}\\
F_{6}&={q+\frac{4pC_{2,1}}{L2}+10\frac{C_{2,1}^{2}}{L2^{2}}}\\
F_{7}&={p+\frac{5C_{2,1}}{L2}}\\
\end{align*}\\
Note: Some of the items above are exhibited as rational functions.  In the table below, these can be used to form polynomial conditions by replacing the rational function by the product of its numerator and denominator.\\
\textbf{Real Root Configurations.}\\
\begin{tabular}{l  p{5cm} p{6cm}}
1.  & 5 distinct real roots & $L3 >0$ $L2 > 0$ and $L1 > 0$ and $D>0$\\
2.  & 3 distinct real roots and 2 complex conjugate roots & $D<0$\\
3.   & 1 real root and 4 distinct complex single roots & $D>0$ and ($L3\leq0$ or $L2\leq0$ or ${L1\leq0}$)\\
4. & 1 double real root and 3 single real roots & $D=0$ and $L1>0$ \\
4.a. & single $<$ double $<$ single $<$ single & $F_2>0$ and $(F_1 \leq 0$ or $F_3 \leq 0)$ \\
4.b. & double $<$ single $<$ single $<$ single & $F_1>0$ and $F_2<0$ and $F_3<0$ \\
4.c. & single $<$ single $<$ single $<$ double & $F_1>0$ and $F_2>0$ and $F_3>0$\\
4.d. & single $<$ single $<$ double $<$ single & $F_2<0$ and $(F_1 \leq 0$ or $F_3 \geq 0)$\\
5.  & 1 double real root and 1 single real root and 2 complex conjugate roots & $D=0$ and $L1<0$\\
5.a. & single $<$ double & $F_2>0$\\
5.b. & double $<$ single & $F_2<0$\\
6.  & 2 real double roots and 1 real single root & $D=0$ and $L1=0$ and $L2\neq0$ and $D2>0$\\
6.a. & single $<$ double $<$ double & $F_4>0$ and $F_5<0$\\
6.b. & double $<$ single $<$ double & $F_4<0$\\
6.c. & double $<$ double $<$ single & $F_4>0$ and $F_5>0$\\
7.  & 2 complex conjugate double roots and 1 single real root & $D=0$ and $L1=0$ and $L2\neq0$ and $D2<0$\\
8.  & 1 triple root and 2 single real roots & $D=0$ and $L1=0$ and $L2>0$ and $D2=0$\\
8.a. & triple $<$ single $<$ single & $F_6>0$ and $F_7<0$\\
8.b. & single $<$ triple $<$ single & $F_6<0$\\
8.c. & single $<$ single $<$ triple & $F_6>0$ and $F_7>0$\\
9.  & 1 triple root and 2 complex conjugate roots & $D=0$ and $L1=0$ and $L2<0$ and $D2=0$\\
\end{tabular}
\newpage
\begin{tabular}{l  p{5cm} p{6cm}}
10.  & 1 quadruple root and 1 single root & $D=0$ and $L1=0$ and $L2=0$ and $L3\neq0$ and $M1=0$\\
10.a. & quadruple $<$ single & $C_4>0$\\
10.b. & single $<$ quadruple & $C_4<0$\\
11.  & 1 triple root and 1 double root & $D=0$ and $L1=0$ and $L2=0$ and $L3\neq0$ and $M1\neq0$\\
11.a. & triple $<$ double & $C_5>0$\\
11.b. & double $<$ triple & $C_5<0$\\
12.  & 1 quintuple root & $D=0$ and $L1=0$ and $L2=0$ and $L3=0$\\
\end{tabular}\bigskip

\textbf{Complex Root Multiplicities.}\\
\begin{tabular}{l  p{5cm} p{6cm}}
1. & 5 distinct roots & $D\neq0$\\
2. & 1 double root and 3 single roots & $D=0$ and $L1\neq0$\\
3. & 2 double roots and 1 single root & $D=0$ and $L1=0$ and $L2\neq0$ and $D2\neq0$\\
4. & 1 triple root and 2 single roots & $D=0$ and $L1=0$ and $L2\neq0$ and $D2=0$\\
5. & 1 quadruple root and 1 single root & $D=0$ and $L1=0$ and $L2=0$ and $M1=0$\\
6. & 1 triple root and 1 double root & $D=0$ and $L1=0$ and $L2=0$ and $M1\neq0$\\
7. & 1 quintuple root & $D=0$ and $L1=0$ and $L2=0$ and $L3=0$\\
\end{tabular}
\newpage

\newpage

%%%%%%%%%%%%%%%%%%%%%%%%%%%%%%%%%%%%%%%%%%%%%%%%%%%%%%%%%%%%%%%%%%%%%%
\section{Bibliography}
%%%%%%%%%%%%%%%%%%%%%%%%%%%%%%%%%%%%%%%%%%%%%%%%%%%%%%%%%%%%%%%%%%%%%%
\noindent [1] Basu, Saugata, Pollock, Richard, and Roy, Marie-Francoise. \emph{Algorithms in Real Algebraic Geometry}. First Edition. Springer-Verlag, Berlin, Heidelberg, 2003.\\
\noindent [2] Benedetti, Riccardo and Risler, Jean-Jacques. \emph{Real Algebraic and Semi-algebraic Sets}. Hermann, Editeurs des Sciences et des Arts, Paris, 1990.\\
\noindent [3] Gonzalez, Eli and Weinberg, David A. ``Root configurations of real univariate cubics and quartics''. arXiv: 1511.07489v2, [math.AC], 8 Jan 2018.\\
\noindent [4] Jacobson, Nathan. \emph{Basic Algebra I}. 2nd ed. W.H. Freeman and Company. 1985.\\
\noindent [5] Liang, S. and Zhang. J. ``A complete discrimination system for polynomials with complex coefficients and its automatic generation.'' \emph{Science in China E} Vol. 42, No. 2, April 1999. p. 113-128.\\
\noindent [6] Sottile, Frank. \emph{Real Solutions to Equations from Geometry}. American Mathematical Society. 2011.\\
\noindent [7] Weinberg, David and Martin, Clyde. ``A note on resultants.'' \emph{Applied Mathematics and Computation} Vol. 24, 1987. p. 303-309.\\
\noindent [8] Yang, L., Hou, X.R., and Zeng, Z.B. ``A complete discrimination system for polynomials.'' \emph{Science in China E} Vol. 39, No. 6, Dec. 1996. p. 628-646.\\
\noindent [9] Yang, Lu. ``Recent advances on determining the number of real roots of paramateric polynomials.'' \emph{Journal of Symbolic Computation} Vol. 28, 1999. p. 225-242.\\

\end{document}